\documentclass[a4paper,12pt]{amsart}
\usepackage{amssymb,amsmath,amsfonts,amscd,euscript}
\usepackage{hyperref}
\usepackage{color}

\newtheorem{theorem}{Theorem}[section]

\theoremstyle{definition}
\newtheorem{definition}[theorem]{Definition}
\makeatother

\usepackage{graphicx}

\begin{document}
	
	\title[PBW degenerations]{PBW degenerations, quiver Grassmannians, and toric varieties}
	\author{Evgeny Feigin}
	\address{HSE University\\
		Faculty of Mathematics\\
		Usacheva 6\\Moscow 119048\\Russia\newline
		{\it and }\newline
		Skolkovo Institute of Science and Technology\\ 
		Center for Advanced Studies\\
		Bolshoy Boulevard 30, bld. 1\\
		Moscow 121205\\
		Russia}
	\email{evgfeig@gmail.com}

\dedicatory{Dedicated to the memory of Ernest Borisovich Vinberg}

	\begin{abstract}
		We present a review on the recently discovered link between the Lie theory,
		the theory of quiver Grassmannians, and various degenerations of flag varieties.
		Our starting point is the induced Poincar\'{e}--Birkhoff--Witt filtration
		on the highest weight representations and the corresponding PBW degenerate
		flag varieties.
	\end{abstract}
	
	\maketitle
	
	\section{Introduction}
	\label{81.sec1}
	
	The celebrated Poincar\'{e}--Birkhoff--Witt theorem claims that there exists
	a filtration on the universal enveloping of a Lie algebra such that the
	associated graded algebra is isomorphic to the symmetric algebra. The PBW
	filtration on the universal enveloping algebra of a nilpotent subalgebra of a simple
	Lie algebra induces a filtration on the representation space of a highest
	weight module. The natural problem is to study this filtration and the
	corresponding graded space. Quite unexpectedly, the problem turned out
	to be related to numerous repre\-sen\-tation-theoretic, algebro-geometric,
	and combinatorial questions. Our goal is to give an overview of the whole
	story and to describe various links between different parts of the picture.
	The main objects of study are monomial bases, convex polytopes, flag and
	Schubert varieties, their degenerations, quiver Grassmannians, and
	toric varieties.
	
	The paper is organized as follows. In~Section~\ref{RTa} we collect re\-pre\-sen\-tation-theoretic
	results of algebraic nature. Section~\ref{RTg} is devoted to the geometric
	representation theory. In~Section~\ref{TC} we discuss combinatorics emerging
	from the cellular decomposition of the PBW degenerate flag varieties. In~Section~\ref{QG} we describe the link between the Lie theory and the theory
	of quiver Grassmannians. Finally, Section~\ref{toric} treats toric degenerations.
	
	Throughout the paper we work over the field of complex numbers.
	
	\section{Representation theory: algebra}
	\label{RTa}
	
	Let ${\mathfrak{g}}$ be a simple Lie algebra with the set $R^{+}$ of positive
	roots. Let $\alpha _{i}$, $\omega _{i}$, $i=1,\dots ,n-1$ be the simple
	roots and the fundamental weights. Let
	${\mathfrak{g}}={\mathfrak{n}}^{+}\oplus {\mathfrak{h}}\oplus {
		\mathfrak{n}}^{-}$ be the Cartan decomposition. 
	For $\alpha \in R^{+}$, we denote
	by $e_{\alpha }\in {\mathfrak{n}}^{+}$ and
	$f_{\alpha }\in {\mathfrak{n}}^{-}$ the corresponding Chevalley generators.
	We denote by $P^{+}$ the set of dominant integral weights.
	
	Consider the PBW filtration on the universal enveloping algebra
	${\mathrm{U}}({\mathfrak{n}}^{-})$:
	\begin{equation*}
		{\mathrm{U}}({\mathfrak{n}}^{-})_{s}=
		\operatorname{span}\bigl\{x_{1}\cdots x_{l}:
		x_{i}\in {\mathfrak{n}}^{-}, l \le s \bigr\},
	\end{equation*}
	for example, ${\mathrm{U}}({\mathfrak{n}}^{-})_{0}={\mathbb{C}}\cdot 1$.
	
	For a dominant integral weight
	$\lambda =m_{1}\omega _{1} + \dots + m_{n-1}\omega _{n-1}$, let
	$V_{\lambda }$ be the corresponding irreducible highest weight
	${\mathfrak{g}}$-module with a highest weight vector $v_{\lambda }$. Since
	$V_{\lambda }={\mathrm{U}}({\mathfrak{n}}^{-})v_{\lambda }$, we have an increasing
	filtration $(V_{\lambda })_{s}$ on $V_{\lambda }$,
	\begin{equation*}
		(V_{\lambda })_{s}={\mathrm{U}}({\mathfrak{n}}^{-})_{s} v_{\lambda }.
	\end{equation*}
	We call this filtration the PBW filtration and study the associated graded
	space $V^{a}_{\lambda }={\mathrm{g}r} V_{\lambda }$.
	
	Let us consider an example for fundamental weights in type $A$. Let
	$V_{\omega _{1}}$ be the vector representation of
	$\mathfrak{sl}_{n}$ with a basis $w_{1},\dots ,w_{n}$ and consider
	$V_{\omega _{k}}\simeq \Lambda ^{k} V_{\omega _{1}}$ for
	$k=1,\dots ,n-1$. Then $(V_{\omega _{k}})_{s}$ is spanned by the wedge
	products $w_{i_{1}}\wedge \dots \wedge w_{i_{k}}$ such that the number
	of indices $a$ with $i_{a}>k$ is at most $s$.
	
	The following holds true \cite{FFL1}:
	\begin{enumerate}
		\item The action of ${\mathrm{U}}({\mathfrak{n}}^{-})$ on $V_{\lambda }$ induces
		the structure of an $S({\mathfrak{n}}^{-})$-module on
		$V_{\lambda }^{a}$ and $V^{a}_{\lambda }=S({\mathfrak{n}}^{-})v_{\lambda }$.
		\item The action of ${\mathrm{U}}({\mathfrak{n}}^{+})$ on $V_{\lambda }$ induces
		the structure of a $U({\mathfrak{n}}^{+})$-module on $V_{\lambda }^{a}$.
	\end{enumerate}
	
	Our aims are to describe $V_{\lambda }^{a}$ as an
	$S({\mathfrak{n}}^{-})$-module and to find a basis of
	$V_{\lambda }^{a}$. We present the answer in type $A$. For similar results
	in other types, see \cite{BK,FFL2,FFL3,G1,G2,M2}.
	
	The positive roots in type $A_{n-1}$ are of the form
	$\alpha _{i,j}=\alpha _{i}+\dots +\alpha _{j}$ with
	$1\le i\le j\le n-1$. Recall that a Dyck path is a sequence
	${\mathbf{p}}=(\beta (0), \beta (1),\dots , \beta (k))$ of positive roots
	of $\mathfrak{sl}_{n}$ satisfying the following conditions: if
	$k=0$, then ${\mathbf{p}}$ is of the form ${\mathbf{p}}=(\alpha _{i})$ for
	some simple root $\alpha _{i}$, and if $k\ge 1$, then the first and last
	elements are simple roots, and if $\beta (s)=\alpha _{p,q}$, then
	$\beta (s+1)=\alpha _{p,q+1}$ or $\beta (s+1)=\alpha _{p+1,q}$.
	
	Here is an example of a path for $\mathfrak{sl}_{6}$:
	\begin{equation*}
		(\alpha _{2},\alpha _{2}+\alpha _{3},
		\alpha _{2}+\alpha _{3}+\alpha _{4}, \alpha
		_{3}+\alpha _{4},\alpha _{4}, \alpha
		_{4}+\alpha _{5},\alpha _{5}).
	\end{equation*}
	
	For a multiexponent ${\mathbf{s}}=\{s_{\beta }\}_{\beta >0}$,
	$s_{\beta }\in {\mathbb{Z}}_{\ge 0}$, let
	$f^{\mathbf{s}}=\prod_{\beta \in R^{+}} f_{\beta }^{s_{\beta }}\in S({
		\mathfrak{n}}^{-})$. For an integral dominant $\mathfrak{sl}_{n}$-weight
	$\lambda =\sum_{i=1}^{n-1} m_{i}\omega _{i}$, let $S(\lambda )$ be the
	set of all multiexponents
	${\mathbf{s}}=(s_{\beta })_{\beta \in R^{+}}\in {\mathbb{Z}}_{\ge 0}^{R^{+}}$
	such that for all Dyck paths
	${\mathbf{p}}=(\beta (0),\dots , \beta (k))$,
	%
	\begin{align}
		\label{upperbound} s_{\beta (0)}+s_{\beta (1)}+\dots + s_{\beta (k)}
		\le m_{i}+m_{i+1}+ \dots + m_{j},
	\end{align}
	where $\beta (0)=\alpha _{i}$ and $\beta (k)=\alpha _{j}$.
	
	The polytopes in ${\mathbb{R}}_{\ge }^{R^{+}}$ defined by inequalities~\eqref{upperbound} are referred to as the FFLV polytopes. For their combinatorial
	properties and connection to the Gelfand--Tsetlin polytopes \cite{GT},
	see \cite{ABS,FMakh,FFLP,FaFo}. The following theorem holds true
	\cite{FFL1}.
	
	\begin{theorem}%
		\label{basis}
		The vectors $f^{\mathbf{s}}v_{\lambda }$, ${\mathbf{s}}\in S(\lambda )$, form
		a basis of $V^{a}_{\lambda }$. In addition,
		$S(\lambda )+S(\mu )=S(\lambda +\mu )$.
	\end{theorem}
	
	We note that Theorem~\ref{basis} implies that the elements
	$f^{\mathbf{s}}v_{\lambda }$, ${\mathbf{s}}\in S(\lambda )$ form a basis of
	the classical representation $V_{\lambda }$ provided an order of factors
	is fixed in each monomial $f^{\mathbf{s}}$ (see \cite{Vin}).
	
	Let us describe the Lie algebra ${\mathfrak{g}}^{a}$ acting on
	$V_{\lambda }^{a}$. As a vector space, ${\mathfrak{g}}^{a}$ is isomorphic
	to ${\mathfrak{g}}$. The Borel
	${\mathfrak{b}}\subset {\mathfrak{g}}^{a}$ is a subalgebra, the nilpotent
	subalgebra ${\mathfrak{n}}^{-}\subset {\mathfrak{g}}^{a}$ is an abelian ideal,
	and ${\mathfrak{b}}$ acts on the space ${\mathfrak{n}}^{-}$ as on the quotient
	${\mathfrak{g}}/{\mathfrak{b}}$. Then for any $\lambda \in P^{+}$ the structure
	of the ${\mathfrak{g}}$-module on $V_{\lambda }$ induces the structures of
	${\mathfrak{g}}^{a}$ module on $V_{\lambda }^{a}$.
	
	Note that $V^{a}_{\lambda }=S({\mathfrak{n}}^{-})v_{\lambda }$ is a cyclic
	$S({\mathfrak{n}}^{-})$-module, so we can write
	$V_{\lambda }^{a}\simeq S({\mathfrak{n}}^{-})/I(\lambda )$, for some ideal
	$I(\lambda )\subset S({\mathfrak{n}}^{-})$.
	
	The following theorem holds in types $A$ and $C$ \cite{FFL1}:
	%
	\begin{theorem}
		$I(\lambda )=S({\mathfrak{n}}^{-})  (\mathrm{U}({\mathfrak{n}}^{+})
		\circ \operatorname{span}\{f_{\alpha }^{(\lambda ,\alpha^\vee )+1}, \alpha \in R^{+}
		\}  )$.
	\end{theorem}
	This theorem should be understood as a commutative analogue of the well-known
	description of $V_{\lambda }$ as the quotient
	\begin{equation*}
		V_{\lambda }\simeq \mathrm{U}({\mathfrak{n}}^{-})/
		\bigl\langle f_{\alpha }^{(
			\lambda ,\alpha^\vee )+1}, \alpha \in R^{+}
		\bigr\rangle
	\end{equation*}
	(see, for example, \cite{FH,J}).
	
	The proof of the theorems above is based on the following claim available
	in types~$A$ and $C$ \cite{FFL1,FFL2,FFL3}.
	
	\begin{theorem}%
		\label{emb}
		Let $\lambda ,\mu \in P^{+}$. Then
		\begin{equation*}
			V_{\lambda +\mu }^{a}\simeq {\mathrm{U}}({
				\mathfrak{g}}^{a}) (v_{\lambda
			}\otimes v_{\mu })
			\hookrightarrow V_{\lambda }^{a}\otimes V_{\mu }^{a}
		\end{equation*}
		as ${\mathfrak{g}}^{a}$-modules.
	\end{theorem}
	
	The algebraic and representation theoretic properties of the PBW filtration
	and the ${\mathfrak{g}}^{a}$ action in more general settings are considered
	in \cite{BD,BFF,BK,CF,Fe1,Fe2,Fe6,FK,FFR,FMake1,G1,G2,M2,PY1,PY2}.
	
	\section{Representation theory: geometry}
	\label{RTg}
	
	Let $G$ be a simple simply-connected Lie group with the Lie algebra
	${\mathfrak{g}}$. Let $B\subset G$ be a Borel subgroup with the Lie algebra
	${\mathfrak{b}}$. Each space $V_{\lambda }$, $\lambda \in P^{+}$ is equipped
	with the natural structure of a $G$-module. Therefore $G$ acts on the projectivization
	${\mathbb{P}}(V_{\lambda })$. The (generalized) flag variety
	$\EuScript{F}_{\lambda }\hookrightarrow {\mathbb{P}}(V_{\lambda })$ is defined
	as the $G$-orbit of the line ${\mathbb{C}}v_{\lambda }$ (see
	\cite{Fu,Kum}). Each variety $\EuScript{F}_{\lambda }$ is isomorphic to
	the quotient of $G$ by the parabolic subgroup leaving the point
	${\mathbb{C}}v_{\lambda }\in {\mathbb{P}}(V_{\lambda })$ invariant. In particular,
	for $G=\operatorname{SL}_{n}$ and a fundamental weight $\lambda =\omega _{k}$ the flag
	variety $\EuScript{F}_{\lambda }$ is isomorphic to the Grassmannian
	$\operatorname{Gr}_{k}(n)$. For a regular weight $\lambda $, the flag variety
	$\EuScript{F}_{\lambda }$ sits inside
	$\prod_{k=1}^{n-1} \operatorname{Gr}_{k}(n)$ and consists of chains of embedded
	subspaces. In what follows, we mostly consider the case $G=\operatorname{SL}_{n}$ and
	regular $\lambda $, the general type $A$ case can be treated similarly
	(see \cite{Fe3,Fe4,Fe6,FeFi}). We denote the complete type $A_{n-1}$ flag
	variety by $\EuScript{F}_{n}$ (it is known to be independent of a regular
	weight $\lambda $). The variety $\EuScript{F}_{n}$ admits Pl\"{u}cker embedding
	into the product of projective spaces
	$\prod_{k=1}^{n-1} {\mathbb{P}}(\Lambda ^{k}({\mathbb{C}}^{n}))$. The homogeneous
	coordinate ring (also known as the Pl\"{u}cker algebra) is a quotient of
	the polynomial algebra in Pl\"{u}cker variables $X_{I}$,
	$I\subset [n]$ by the quadratic Pl\"{u}cker ideal.
	
	Recall the Lie algebra ${\mathfrak{g}}^{a}$ acting on
	$V_{\lambda }^{a}$. We now describe the corresponding Lie group
	$G^{a}$. Let $M=\dim {\mathfrak{n}}$ and let ${\mathbb{G}}_{a}$ be the additive
	group of the field ${\mathbb{C}}$. The Lie group $G^{a}$ is a semidirect
	product ${\mathbb{G}}_{a}^{M}\rtimes B$ of the normal subgroup
	${\mathbb{G}}_{a}^{M}$ and the Borel subgroup~$B$. The action by conjugation
	of $B$ on ${\mathbb{G}}_{a}^{M}$ is induced from the $B$-action on
	$({\mathfrak{n}}^{-})^{a}\simeq {\mathfrak{g}}/{\mathfrak{b}}$.
	
	We now define the degenerate flag varieties
	$\EuScript{F}_{\lambda }^{a}$ \cite{Fe3}. Let
	$[v_{\lambda }]\in {\mathbb{P}}(V_{\lambda }^{a})$ be the line~${\mathbb{C}}v_{\lambda }$.
	%
	\begin{definition}
		The variety
		$\EuScript{F}_{\lambda }^{a}\hookrightarrow {\mathbb{P}}(V_{\lambda }^{a})$
		is the closure of the $G^{a}$-orbit of $[v_{\lambda }]$,
		\begin{equation*}
			\EuScript{F}_{\lambda }^{a}=\overline{G^{a}
				[v_{\lambda }]}= \overline{{\mathbb{G}}_{a}^{M}
				[v_{\lambda }]}\hookrightarrow {\mathbb{P}}(V_{\lambda }^{a}).
		\end{equation*}
	\end{definition}
	
	We note that the orbit
	$G [v_{\lambda }]\hookrightarrow {\mathbb{P}}(V_{\lambda })$ coincides with
	its closure, but the orbit $G^{a} [v_{\lambda }]$ does not; in fact,
	$\EuScript{F}_{\lambda }^{a}$ is the so-called ${\mathbb{G}}_{a}^{M}$-variety,
	see \cite{HT,Ar,AS}. Theorem~\ref{emb} implies that in types $A$ and
	$C$ the varieties $\EuScript{F}_{\lambda }^{a}$ depend only on the regularity
	class of $\lambda $, i.e., $\EuScript{F}_{\lambda }^{a}$~is isomorphic to
	$\EuScript{F}_{\mu }^{a}$ if and only if the sets of fundamental weights
	showing up in $\lambda $ and $\mu $ coincide (see \cite{M1} for the study
	of a similar question for Schubert varieties).
	
	In types $A$ and $C$, we have rather explicit description of the degenerate
	flag varieties \cite{Fe4,FFLi}. In particular, for
	${\mathfrak{g}}=\mathfrak{sl}_{n}$ one has
	$\EuScript{F}_{\omega _{k}}^{a}\simeq \operatorname{Gr}_{k}(n)$. To describe
	the PBW degenerate flag varieties in type $A$, we introduce the following
	notation: let $W$ be an $n$-dimensional vector space with a basis
	$w_{1},\dots ,w_{n}$. Let us denote by $\operatorname{pr}_{k}: W\to W$ the projection
	along $w_{k}$. We denote the regular PBW degenerate flag variety by
	$\EuScript{F}^{\,a}_{n}$. The following theorem holds \cite{Fe3,Fe4} (we
	use the shorthand notation $[n]=\{1,\dots ,n\}$).
	
	\begin{theorem}%
		\label{explicit}
		One has
		\begin{equation*}
			\EuScript{F}_{n}^{a}\simeq \bigl\{(V_{1},
			\dots ,V_{n-1}): V_{k}\in \operatorname{Gr}_{k}(W),
			k\in [n]; \operatorname{pr}_{k+1}V_{k}\subset
			V_{k+1}, k\in [n-1] \bigr\}.
		\end{equation*}
	\end{theorem}
	
	Using this description, one proves the following theorem
	\cite{CL,CLL,CFR4} (see also \cite{LS}).
	%
	\begin{theorem}
		The variety $\EuScript{F}^{\,a}_{n}$ is isomorphic to a Schubert variety
		for the group $\operatorname{SL}_{2n-1}$.
	\end{theorem}
	
	The symplectic PBW degenerations are described in \cite{FFLi} (see also
	\cite{BF2}).
	
	For a partition
	$\lambda =(\lambda _{1}\ge \dots \ge \lambda _{n-1}\ge 0)$, we denote by
	$Y_{\lambda }$ the corresponding Young diagram. Recall that the classical
	$\operatorname{SL}_{n}$ flag variety admits an embedding to the product of Grassmannians.
	The corresponding homogeneous coordinate ring (the Pl\"{u}cker algebra) is
	generated by the Pl\"{u}cker variables $X_{I}$, $I\subset [n]$ and is known
	to be isomorphic to the direct sum
	$\bigoplus_{\lambda \in P^{+}} V_{\lambda }^{*}$ (see \cite{Fu}). There
	is a one-to-one bijection between the Pl\"{u}cker variables and columns filled
	with numbers from $[n]$ (the numbers increase from top to bottom).
	Then the semistandard Young tableaux provide a basis of the homogeneous
	coordinate ring of $\operatorname{SL}_{n}/B$ (one takes the product of Pl\"{u}cker variables,
	corresponding to the columns of a tableau). Similar result holds true in
	the PBW degenerate situation.
	
	We denote by $\mu _{j}$ the length of the $j$th column of a diagram.
	
	\begin{definition}%
		\label{t}
		A semistandard PBW-tableau of shape $\lambda $ is a filling
		$T_{i,j}$ of the Young diagram $Y_{\lambda }$ with numbers
		$1,\dots ,n$. The number $T_{i,j}\in \{1,\dots , n\}$ is attached to the
		box $(i,j)$. The filling $T_{i,j}$ has to satisfy the following properties:
		\begin{enumerate}
			\item if $T_{i,j}\le \mu _{j}$, then $T_{i,j}=i$;
			\item if $i_{1}<i_{2}$ and $T_{i_{1},j}\ne i_{1}$, then
			$T_{i_{1},j} > T_{i_{2},j}$;
			\item for any $j>1$ and any $i$, there exists $i_{1}\ge i$ such that
			$T_{i_{1},j-1}\ge T_{i,j}$.
		\end{enumerate}
	\end{definition}
	
	One can show that the number of shape $\lambda $ semistandard PBW-tableaux
	is equal to $\dim V_{\lambda }$. Moreover, the following theorem holds
	\cite{Fe3} (see also \cite{Ha,FFL4}).
	
	\begin{theorem}
		The homogeneous coordinate ring of $\EuScript{F}^{\,a}_{n}$ (also known
		as the PBW degenerate Pl\"{u}cker algebra) is isomorphic to the direct sum
		of dual PBW degenerate modules $(V_{\lambda }^{a})^{*}$,
		$\lambda \in P^{+}$. The ideal of relations is quadratic and is generated by degenerate Pl\"{u}cker relations. The PBW semistandard tableaux parametrize
		a basis in the coordinate ring.
	\end{theorem}
	
	Certain infinite-dimensional analogues of the results described above are
	obtained in \cite{FeFiR,P}. However, this direction has not been seriously
	pursued so far.
	
	\section{Topology and combinatorics}
	\label{TC}
	
	In this section we describe a cellular decomposition of the type $A$ complete
	PBW degenerate flag varieties $\EuScript{F}^{\,a}_{n}$ (see
	\cite{Fe4,BF2,FFLi} for a more general picture).
	
	Let us fix an $n$-dimensional vector space $W$ with a basis
	$w_{1},\dots ,w_{n}$. Let ${{\mathbf{I}}}=(I_{1},\dots ,I_{n-1})$ be a collection
	of subsets of the set $[n]$ such that $|I_{k}|=k$. We denote by
	$p_{{\mathbf{I}}}\in \prod_{k=1}^{n-1} \operatorname{Gr}_{k}(W)$ a point in
	the product of Grassmann varieties such that the $k$th component is equal
	to the linear span of $w_{i}$ with $i\in I_{k}$. Theorem~\ref{explicit} implies that
	$p_{{\mathbf{I}}}\in \EuScript{F}^{\,a}_{n}$ if and only if
	%
	\begin{align}
		\label{Genocchi} I_{k}\subset I_{k+1}\cup \{k+1\}
		\quad \text{for all } k=1,\dots ,n-2.
	\end{align}
	
	The following theorem is proved in \cite{Fe4}.
	%
	\begin{theorem}
		The $G^{a}$ orbits of the points $p_{\mathbf{I}}$ provide a cellular decomposition
		of $\EuScript{F}^{\,a}_{n}$.
	\end{theorem}
	
	A natural problem is to compute the Euler characteristic and Poincar\'{e} polynomial
	of $\EuScript{F}^{\,a}_{n}$. The answer is given in terms of the normalized
	median Genocchi numbers and the Dellac configurations.
	
	The normalized median Genocchi numbers $h_{n}$, $n=0,1,2,\dots $ form a
	sequence which starts with $1,1,2,7,38,295$ \cite{OEIS}. The earliest definition
	was given by Dellac in \cite{De} (see also
	\cite{B,BF1,Ba,Du,DR,DV,Fe5,HZ,Kr,Vie,ZZ}). Consider a rectangle with
	$n$ columns and $2n$ rows. It contains $n\times 2n$ boxes labeled by pairs
	$(l,j)$, where $l=1,\dots ,n$ is the number of a column and
	$j=1,\dots ,2n$ is the number of a row. A Dellac configuration $D$ is a
	subset of boxes, subject to the following conditions:
	\begin{enumerate}
		\item each column contains exactly two boxes from $D$,
		\item each row contains exactly one box from $D$,
		\item if the $(l,j)$th box is in $D$, then $l\le j\le n+l$.
	\end{enumerate}
	Let $DC_{n}$ be the set of such configurations. Then the number of elements
	in $DC_{n}$ is equal to~$h_{n}$.
	
	We list all Dellac's configurations for $n=3$.
	
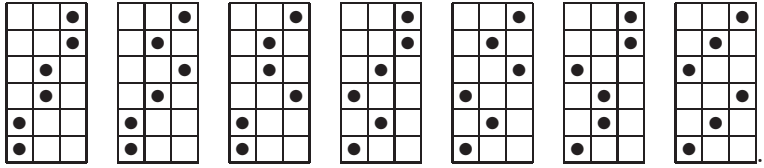
\begin{equation*}
		\begin{picture}(30,60)
			\put(0,0){\line(1,0){30}}
			\put(0,10){\line(1,0){30}}
			\put(0,20){\line(1,0){30}}
			\put(0,30){\line(1,0){30}}
			\put(0,40){\line(1,0){30}}
			\put(0,50){\line(1,0){30}}
			\put(0,60){\line(1,0){30}}
			
			\put(0,0){\line(0,1){60}}
			\put(10,0){\line(0,1){60}}
			\put(20,0){\line(0,1){60}}
			\put(30,0){\line(0,1){60}}
			
			\put(2,2){$\bullet $}
			\put(2,12){$\bullet $}
			\put(12,22){$\bullet $}
			\put(12,32){$\bullet $}
			\put(22,42){$\bullet $}
			\put(22,52){$\bullet $}
		\end{picture} %
		\quad \begin{picture}(30,60)
			\put(0,0){\line(1,0){30}}
			\put(0,10){\line(1,0){30}}
			\put(0,20){\line(1,0){30}}
			\put(0,30){\line(1,0){30}}
			\put(0,40){\line(1,0){30}}
			\put(0,50){\line(1,0){30}}
			\put(0,60){\line(1,0){30}}
			
			\put(0,0){\line(0,1){60}}
			\put(10,0){\line(0,1){60}}
			\put(20,0){\line(0,1){60}}
			\put(30,0){\line(0,1){60}}
			
			\put(2,2){$\bullet $}
			\put(2,12){$\bullet $}
			\put(12,22){$\bullet $}
			\put(12,42){$\bullet $}
			\put(22,32){$\bullet $}
			\put(22,52){$\bullet $}
		\end{picture} %
		\quad \begin{picture}(30,60)
			\put(0,0){\line(1,0){30}}
			\put(0,10){\line(1,0){30}}
			\put(0,20){\line(1,0){30}}
			\put(0,30){\line(1,0){30}}
			\put(0,40){\line(1,0){30}}
			\put(0,50){\line(1,0){30}}
			\put(0,60){\line(1,0){30}}
			
			\put(0,0){\line(0,1){60}}
			\put(10,0){\line(0,1){60}}
			\put(20,0){\line(0,1){60}}
			\put(30,0){\line(0,1){60}}
			
			\put(2,2){$\bullet $}
			\put(2,12){$\bullet $}
			\put(12,32){$\bullet $}
			\put(12,42){$\bullet $}
			\put(22,22){$\bullet $}
			\put(22,52){$\bullet $}
		\end{picture}
		\quad \begin{picture}(30,60)
			\put(0,0){\line(1,0){30}}
			\put(0,10){\line(1,0){30}}
			\put(0,20){\line(1,0){30}}
			\put(0,30){\line(1,0){30}}
			\put(0,40){\line(1,0){30}}
			\put(0,50){\line(1,0){30}}
			\put(0,60){\line(1,0){30}}
			
			\put(0,0){\line(0,1){60}}
			\put(10,0){\line(0,1){60}}
			\put(20,0){\line(0,1){60}}
			\put(30,0){\line(0,1){60}}
			
			\put(2,2){$\bullet $}
			\put(2,22){$\bullet $}
			\put(12,12){$\bullet $}
			\put(12,32){$\bullet $}
			\put(22,42){$\bullet $}
			\put(22,52){$\bullet $}
		\end{picture}%
		\quad \begin{picture}(30,60)
			\put(0,0){\line(1,0){30}}
			\put(0,10){\line(1,0){30}}
			\put(0,20){\line(1,0){30}}
			\put(0,30){\line(1,0){30}}
			\put(0,40){\line(1,0){30}}
			\put(0,50){\line(1,0){30}}
			\put(0,60){\line(1,0){30}}
			
			\put(0,0){\line(0,1){60}}
			\put(10,0){\line(0,1){60}}
			\put(20,0){\line(0,1){60}}
			\put(30,0){\line(0,1){60}}
			
			\put(2,2){$\bullet $}
			\put(2,22){$\bullet $}
			\put(12,12){$\bullet $}
			\put(12,42){$\bullet $}
			\put(22,32){$\bullet $}
			\put(22,52){$\bullet $}
		\end{picture}%
		\quad \begin{picture}(30,60)
			\put(0,0){\line(1,0){30}}
			\put(0,10){\line(1,0){30}}
			\put(0,20){\line(1,0){30}}
			\put(0,30){\line(1,0){30}}
			\put(0,40){\line(1,0){30}}
			\put(0,50){\line(1,0){30}}
			\put(0,60){\line(1,0){30}}
			
			\put(0,0){\line(0,1){60}}
			\put(10,0){\line(0,1){60}}
			\put(20,0){\line(0,1){60}}
			\put(30,0){\line(0,1){60}}
			
			\put(2,2){$\bullet $}
			\put(2,32){$\bullet $}
			\put(12,12){$\bullet $}
			\put(12,22){$\bullet $}
			\put(22,42){$\bullet $}
			\put(22,52){$\bullet $}
		\end{picture}%
		\quad \begin{picture}(30,60)
			\put(0,0){\line(1,0){30}}
			\put(0,10){\line(1,0){30}}
			\put(0,20){\line(1,0){30}}
			\put(0,30){\line(1,0){30}}
			\put(0,40){\line(1,0){30}}
			\put(0,50){\line(1,0){30}}
			\put(0,60){\line(1,0){30}}
			
			\put(0,0){\line(0,1){60}}
			\put(10,0){\line(0,1){60}}
			\put(20,0){\line(0,1){60}}
			\put(30,0){\line(0,1){60}}
			
			\put(2,2){$\bullet $}
			\put(2,32){$\bullet $}
			\put(12,12){$\bullet $}
			\put(12,42){$\bullet $}
			\put(22,22){$\bullet $}
			\put(22,52){$\bullet $}
		\end{picture}%
		.
	\end{equation*}
	
	The importance of the median Genocchi numbers comes from the following
	theorem~\cite{Fe4}.
	
	\begin{theorem}
		The number of collections ${{\mathbf{I}}}$ subject to conditions
		\textup{\eqref{Genocchi}} is equal to the normalized median Genocchi number
		$h_{n}$. The Euler characteristic of $\EuScript{F}^{\,a}_{n}$ is equal to
		$h_{n}$.
	\end{theorem}
	
	An explicit formula for the numbers $h_{n}$ is available (see
	\cite{CFR1}), namely
	%
	\begin{align}
		h_{n}=\sum_{f_{0},\dots ,f_{n}\ge 0} \prod
		_{k=1}^{n} \binom{1+f_{k-1}}{f_{k}} \prod
		_{k=0}^{n-1} \binom{1+f_{k+1}}{f_{k}}
	\end{align}
	with $f_{0}=f_{n}=0$.
	
	In order to compute the Poincar\'{e} polynomial of
	$\EuScript{F}^{\,a}_{n}$, we define a length $l(D)$ of a Dellac configuration
	$D$ as the number of pairs $(l_{1},j_{1})$, $(l_{2},j_{2})$ such that the
	boxes $(l_{1},j_{1})$ and $(l_{2},j_{2})$ are both in $D$ and
	$l_{1}<l_{2}$, $j_{1}>j_{2}$. This definition resembles the definition
	of the length of a permutation. We note that in the classical case the
	complex dimension of the cell attached to a permutation $\sigma $ in a
	flag variety is equal to the number of pairs $j_{1}<j_{2}$ such that
	$\sigma (j_{1})>\sigma (j_{2})$, which equals to the length of
	$\sigma $. One has \cite{Fe4}:
	
	\begin{theorem}
		The complex dimension of the cell in $\EuScript{F}^{\,a}_{n}$ containing
		a point $p_{{\mathbf{I}}}$ is equal to $l(D)$. Thus the Poincar\'{e} polynomial
		$h_{n}(q)=P_{\EuScript{F}^{\,a}_{n}}(q)$ is given by
		$h_{n}(q)=\sum_{D\in DC_{n}} q^{l(D)}$.
	\end{theorem}
	The first four polynomials $h_{n}(q)$ are as follows:
	\begin{align*}
		h_{1}(q)&=1, \quad  h_{2}(q)=1+q,
		\\
		h_{3}(q)&=1+ 2q+ 3q^{2}+q^{3},
		\\
		h_{4}(q)&=1+3q+7q^{2}+10q^{3}+10q^{4}+6q^{5}+q^{6}.
	\end{align*}
	
	The following (fermionic type) formula for the polynomials
	$h_{n}(q)$ is obtained in \cite{CFR1} using the geometry of quiver Grassmannians:
	%
	\begin{align}
		h_{n}(q)= \sum_{f_{1},\dots ,f_{n-1}\ge 0}q^{\sum _{k=1}^{n-1} (k-f_{k})(1-f_{k}+f_{k+1})}
		\prod_{k=1}^{n} \binom{1+f_{k-1}}{f_{k}}_{ q}
		\prod_{k=0}^{n-1} \binom{1+f_{k+1}}{f_{k}}_{ q}
	\end{align}
	(we assume $f_{0}=f_{n}=0$). The formula is given in terms of the
	$q$-binomial coefficients
	\begin{equation*}
		\binom{m}{n}_{ q}=\frac{m_{q}!}{n_{q}!(m-n)_{q}!},\quad
		m_{q}!= \prod_{i=1}^{m}
		\frac{1-q^{i}}{1-q}.
	\end{equation*}
	
	\section{Quiver Grassmannians}
	\label{QG}
	
	Theorem~\ref{explicit} provides a link between the PBW degenerate flag
	varieties and quiver Grassmannians. Let $Q$ be a quiver with the set of
	vertices $Q_{0}$ and the set of arrows $Q_{1}$. For two vectors
	${\mathbf{e}},{\mathbf{d}}\in {\mathbb{Z}}^{Q_{0}}$, we denote by
	$\langle {\mathbf{e}},{\mathbf{d}}\rangle $ the value of the Euler from of
	the quiver. For a $Q$ module $M$ and a dimension vector
	${{\mathbf{e}}}\in {\mathbb{Z}}_{\ge 0}^{Q_{0}}$, we denote by
	$\operatorname{Gr}_{\mathbf{e}}(M)$ the quiver Grassmannian consisting of
	${\mathbf{e}}$-dimensional subrepresentations of $M$. For more details on
	the quiver representation theory, see \cite{ASS,CB1,CB2,Schi}. The general
	theory of quiver Grassmannians can be found in \cite{CI} (see also
	\cite{AdF1,AdF2,CR,CEFR,Hubery,LW,REveryProj,RingEveryProj}).
	
	Now let $Q$ be an equioriented type $A_{n-1}$ quiver. We label the vertices
	by the numbers from $1$ to $n$. Then the set $Q_{1}$ consists of arrows
	$i\to i+1$, $i\in [n-1]$. The indecomposable representations of $Q$ are
	labeled by pairs $1\le i\le j\le n$; the representation $U_{i,j}$ is supported
	on vertices from $i$ to $j$ and is one-dimensional at each vertex. The
	projective indecomposable representations are given by
	$P_{k}=U_{k,n}$ and the injective indecomposables are
	$I_{k}=U_{1,k}$. In particular, the path algebra $A$ of $Q$ is isomorphic
	to the direct sum $\bigoplus_{k=1}^{n-1} P_{k}$ and the dual
	$A^{*}$ is the direct sum $\bigoplus_{k=1}^{n-1} I_{k}$ of all indecomposable
	injectives.
	
	By the very definition, the classical complete flag variety
	$\operatorname{SL}_{n}/B$ is isomorphic to the quiver Grassmannian
	$\mathrm{Gr}_{\dim A} (P_{1}^{\oplus n})$. The following observation was
	made in \cite{CFR1}:
	%
	\begin{align}
		\label{df=qg} \EuScript{F}^{\,a}_{n}\simeq
		\mathrm{Gr}_{\dim A} (A\oplus A^{*}).
	\end{align}
	The realization~\eqref{df=qg} provides additional tools for the study of
	algebro-geometric and combinatorial properties of the degenerate flag varieties
	(see \cite{CFR1,CFR2,CFR3,CFR4}). In particular, one recovers and generalizes
	\cite{CFR3,CFR5} the Bott--Samelson type construction for the resolution
	of singularities of $\EuScript{F}^{\,a}_{n}$ \cite{FeFi} (see also
	\cite{KS,Sche} for further generalizations). The resolution is constructed
	as a quiver Grassmannian for a larger quiver attached to $Q$.
	
	Since the degenerate flag varieties have many nice properties, it is natural
	to study the quiver Grassmannians $\mathrm{Gr}_{\dim P} (P\oplus I)$ for
	arbitrary projective representation $P$ and an injective representation
	$I$ and a Dynkin quiver $Q$ (the so-called well-behaved quiver Grassmannians).
	We summarize the main properties of these quiver Grassmannians in the following
	theorem (see \cite{CFR1,CFR2}).
	
	\begin{theorem}
		Let $P$ and $I$ be a projective and an injective representations of a Dynkin
		quiver $Q$. Then the quiver Grassmannian
		$X=\mathrm{Gr}_{\dim P} (P\oplus I)$ has the following properties:
		\begin{enumerate}
			\item $\dim X = \langle \dim P,\dim I\rangle $,
			\item $X$ is irreducible and normal,
			\item $X$ is locally a complete intersection,
			\item there exists an algebraic group
			$G\subset {\mathrm{Aut}}(P\oplus I)$ acting on $X$ with finitely many orbits.
		\end{enumerate}
	\end{theorem}
	
	For a dimension vector ${\mathbf{d}}\in {\mathbb{Z}}_{\ge 0}^{Q_{0}}$, let
	$R_{\mathbf{d}}$ be the variety of $Q$-rep\-resentations of dimension
	${\mathbf{d}}$. The group
	${\mathrm{G}L}_{\mathbf{d}}=\prod_{i\in Q_{0}} {\mathrm{G}L}_{d_{i}}$ acts
	on $R_{\mathbf{d}}$ by base change and the orbits are parameterized by the
	isoclasses of ${\mathbf{d}}$-dimen\-sional representations of $Q$. The closure
	of orbits induces the degeneration order on the set of isoclasses. Fixing
	a dimension vector ${\mathbf{e}}$, we obtain a family
	$\operatorname{Gr}_{\mathbf{e}}({\mathbf{d}})$ of ${\mathbf{e}}$-dimensional quiver
	Grassmannians over the representation space $R_{\mathbf{d}}$ (the so-called
	universal quiver Grassmannian). Let us denote the projection map
	$\operatorname{Gr}_{\mathbf{e}}({\mathbf{d}})\to R_{\mathbf{d}}$ by
	$p_{{\mathbf{e}},{\mathbf{d}}}$.
	
	We are interested in the case when $Q$ is the equioriented type
	$A_{n-1}$ quiver, ${\mathbf{d}}=(n,\dots ,n)$ and
	${\mathbf{e}}=(1,2,\dots ,n-1)$. Then both the classical and the PBW degenerate
	flag varieties are isomorphic to the fibers of
	$p_{{\mathbf{e}},{\mathbf{d}}}$. It is thus natural to ask about the properties
	of the whole family. The ${\mathrm{G}L}_{\mathbf{d}}$ orbits on
	$R_{\mathbf{d}}$ are parametrized by the tuples ${\mathbf{r}}$ of ranks
	$r_{i,j}$ of the compositions of the maps between the $i$th and $j$th vertices.
	We define three rank tuples ${\mathbf{r}}^{0}$, ${\mathbf{r}}^{1}$, and
	${\mathbf{r}}^{2}$ by
	\begin{equation*}
		r^{0}_{i,j}=n+1,\quad  r^{1}_{i,j}=n+1-(j-i),
		\quad  r^{2}_{i,j}=n-(j-i).
	\end{equation*}
	Then the corresponding representations of $Q$ are given by
	$M^{0}=P_{1}^{\oplus n}$, $M^{1}=A\oplus A^{*}$, and
	\begin{equation*}
		M^{2} = \bigoplus_{k=1}^{n-1}
		P_{k}\oplus \bigoplus_{k=1}^{n-2}
		I_{k} \oplus S,
	\end{equation*}
	where $S$ is the direct sum of all simple modules of $Q$. One has
	$\operatorname{SL}_{n}/B\simeq \operatorname{Gr}_{\mathbf{e}}(M^{0})$,
	$\EuScript{F}^{\,a}_{n}\simeq \operatorname{Gr}_{\mathbf{e}}(M^{1})$. In
	\cite{CFFFR1} we prove the following theorem:
	
	\begin{theorem}\phantomsection
		\begin{enumerate}
			\item[(a)] The quiver Grassmannian
			$p_{{\mathbf{e}},{\mathbf{d}}}^{-1}(M^{2})$ is of expected dimension
			$n(n-1)/2$. It is reducible and the number of irreducible components is
			equal to the $n$th Catalan number.
			\item[(b)] The flat irreducible locus of
			$\operatorname{Gr}_{\mathbf{e}}({\mathbf{d}})$ consists of the fibers
			$p_{{\mathbf{e}},{\mathbf{d}}}^{-1}(M)$ such that $M$ degenerates to
			$M^{1}$.
			\item[(c)] The flat locus of $\operatorname{Gr}_{\mathbf{e}}({\mathbf{d}})$ consists
			of the fibers $p_{{\mathbf{e}},{\mathbf{d}}}^{-1}(M)$ such that $M$ degenerates
			to $M^{2}$.
		\end{enumerate}
	\end{theorem}
	
	The case of partial flag varieties is considered in \cite{CFFFR2}.
	
	\section{Toric degenerations}
	\label{toric}
	
	As explained in Section~\ref{QG}, the degeneration of the classical flag
	variety into the PBW degenerate flag variety can be considered within a
	family of quiver Grassmannians over the representation space of the quiver.
	In particular, the study of other degenerations (intermediate and deeper
	ones) leads to the new and interesting results and examples. Yet another
	direction is to make a connection between the PBW degeneration and toric
	degenerations \cite{CLS} of flag varieties (the latter attracted a lot
	of attention in the last two decades, see
	\cite{AB,A,BCKS,BFFHL,Cal,FaFL1,FaFL2,FLP,HK,Lak}). One of the most famous
	examples is a flat degeneration of $\EuScript{F}_{n}$ into the toric variety
	with the Newton polytope being the Gelfand--Tsetlin polytope
	\cite{GT,S,KM,GLa}. We are able to prove the following theorem.
	
	\begin{theorem}%
		\label{toricFFLV}
		The complete flag variety $\EuScript{F}_{n}=\operatorname{SL}_{n}/B$ admits a flat degeneration
		to the toric variety corresponding to a FFLV polytope with regular highest
		weight. This degeneration factors through the PBW degeneration.
	\end{theorem}
	
	The GT and FFLV polytopes are identified with the order and chain polytopes
	of a certain poset (see \cite{ABS,Stanley,FFLP,M4,MY}). Several proofs
	of Theorem~\ref{toricFFLV} are available. Essentially, there are three
	different approaches:
	\begin{enumerate}
		\item via the representation space $V_{\lambda }$,
		\item via the Gr\"{o}bner theory for the Pl\"{u}cker ideal,
		\item via the SAGBI theory for the Pl\"{u}cker algebra.
	\end{enumerate}
	
	The first approach is utilized in \cite{FFL4,FFR,FFFM}. The approach is
	similar to the PBW degeneration construction: instead of attaching degree
	one to each Chevalley generator, one uses a weight system, attaching weight
	$a_{i,j}$ to the generators $f_{\alpha _{i,j}}$ for all positive roots.
	For certain weight systems, one gets a filtration on the universal enveloping
	algebra, which leads to a filtered (and then graded) representation space
	and degenerate flag variety. Here comes the theorem (see
	\cite{FFR,FFL4}).
	
	\begin{theorem}
		Consider the PBW filtration with the weight system
		$a_{i,j}=(j-i+1)(n-j)$. Then
		\begin{enumerate}
			\item in the associated graded space the nonzero monomials in
			$f_{\alpha }$ form a basis,
			\item the associated graded space is acted upon by the symmetric algebra
			$S({\mathfrak{n}}^{-})$ and the degenerate flag variety is a
			${\mathbb{G}}_{a}^{n(n-1)/2}$ variety,
			\item the corresponding degenerate flag variety is toric with the Newton
			polytope being the FFLV polytope.
		\end{enumerate}
	\end{theorem}
	
	Instead of working with the representation space, one may start with the
	algebraic variety $\EuScript{F}_{n}$ from the very beginning. As an intermediate
	step one considers the theory of Newton--Okounkov (NO) bodies
	\cite{O,KK}. The connection between the NO bodies and toric degenerations
	is used in many papers, see, e.g., \cite{A,FL,HK,Ka}. The following holds
	true.
	
	\begin{theorem}
		The toric variety attached to the FFLV polytope can be constructed as a
		Newton--Okounkov body for certain valuations. The valuations are obtained
		via Lie theory \cite{FFL4} or geometrically \cite{FujH,Fuj,Kir1,Kir2}.
	\end{theorem}
	
	Recall the Pl\"{u}cker coordinates $X_{I}$, the quadratic Pl\"{u}cker ideal
	defining the flag variety $\EuScript{F}_{n}$ inside the product of the
	projectivized fundamental representations and the Pl\"{u}cker algebra (the
	quotient by the Pl\"{u}cker ideal). There are two general constructions leading
	to the degenerations of algebraic varieties: Gr\"{o}bner theory for the defining
	ideals \cite{MR,HH} (see also \cite{SpSt,SW,FFFM} for the tropical version)
	and the SAGBI (subalgebra analogues of the Gr\"{o}bner bases for ideals)
	theory \cite{RS} (see also \cite{DEP,HH,Hi}). The former construction works
	with the defining ideals, attaching certain degrees to the variables, and
	the latter deals with the quotient algebras, using certain monomial orders.
	In our setting the following claims hold (see \cite{FFFM,M3}).
	
	\begin{theorem}
		There exists a maximal cone in the Gr\"{o}bner fan of the Pl\"{u}cker ideal
		such that a general point corresponds to the monomial ideal defined by
		the PBW semistandard tableaux. There exists a monomial order on the set
		of Pl\"{u}cker variables such that the monomials in Pl\"{u}cker variables corresponding
		to the PBW semistandard tableaux form a SAGBI basis of the Pl\"{u}cker algebra.
	\end{theorem}
	
	Let us close with the remark that it would be very interesting to construct
	and study toric degenerations for affine flag varieties \cite{Kum} and
	semiinfinite flag varieties \cite{FeFr,FiMi}. The first steps in this direction
	were made in \cite{Sot,SS}. From the representation theory point of view,
	this would lead to new constructions of bases and character formulas for
	the integrable representations of affine algebras and global Weyl and Demazure
	modules for current algebras \cite{DF,DFF,FMake2,Make}.

\section{Acknowledgments}	
		I dedicate this review to the memory of Ernest Borisovich Vinberg, who
		passed away in 2020. I am indebted to him for sharing his ideas on monomial
		bases in irreducible representations of simple Lie algebras. I am grateful
		to Giovanni Cerulli Irelli, Xin Fang, Michael Finkelberg, Ghislain Fourier, Peter Littelmann, Igor Makhlin
		and Markus Reineke for fruitful collaboration.
			This work was partially funded within the framework of the HSE University
		Basic Research Program.
	
	
	

\end{document}